\documentclass{eptcs}

\usepackage{iftex}

\ifpdf
  \usepackage{underscore}         
  \usepackage[T1]{fontenc}        
\else
  \usepackage{breakurl}           
\fi

\usepackage{spalign}
\usepackage{physics}
\usepackage{amsmath}
\usepackage{mathtools}
\usepackage{amssymb}
\usepackage{stmaryrd}
\usepackage{graphicx}
\usepackage[svgnames]{xcolor}
\usepackage{tikz-cd}
\usepackage{tikz}
\usetikzlibrary{automata, arrows.meta, positioning}
\usepackage{amsthm}
\usepackage{rotating}
\usepackage{adjustbox}
\usepackage{hyperref}
\usepackage[capitalise]{cleveref}
\usepackage{import}
\usepackage{changepage}
\usepackage{ebutf8}

\usetikzlibrary{calc}

\theoremstyle{plain}
\newtheorem{theorem}{Theorem}[section]
\newtheorem{lemma}[theorem]{Lemma}
\newtheorem{proposition}[theorem]{Proposition}
\newtheorem{corollary}[theorem]{Corollary}

\theoremstyle{definition}
\newtheorem{definition}[theorem]{Definition}
\newtheorem{example}[theorem]{Example}

\newtheorem{construction}[theorem]{Construction}

\theoremstyle{remark}

\usepackage{xcolor}

\newcommand\imCMsym[4][\mathord]{%
  \DeclareFontFamily{U} {#2}{}
  \DeclareFontShape{U}{#2}{m}{n}{
    <-6> #25
    <6-7> #26
    <7-8> #27
    <8-9> #28
    <9-10> #29
    <10-12> #210
    <12-> #212}{}
  \DeclareSymbolFont{CM#2} {U} {#2}{m}{n}
  \DeclareMathSymbol{#4}{#1}{CM#2}{#3}
}

\imCMsym{cmmi}{124}{\CMjmath}
\imCMsym[\mathop]{cmsy}{113}{\CMamalg}

\tikzset{
rot/.style={anchor=north, rotate=90, inner sep=0pt}
}

\usepackage{xspace}
\newcommand{\Fin}{\textbf{Fin}\xspace}
\DeclareMathOperator{\Set}{\textbf{Set}}
\newcommand*{\Rel}{\textbf{Rel}\xspace}
\newcommand*{\SpanSet}{\textbf{Span(Set)}\xspace}
\newcommand*{\SpanFin}{\textbf{Span(Fin)}\xspace}
\DeclareMathOperator{\Id}{Id}
\renewcommand{\Im}{\mathrm{Im}}
\DeclareMathOperator{\word}{word}

\definecolor{myGreen}{RGB}{6, 150, 10}

\title{Fibrational Perspectives on Determinization of Finite-State Automata}
\author{Thea Li
\institute{Université Paris Cité \\ Paris, France}
\email{thea.li@etu.u-paris.fr}}
\date{}
\def\titlerunning{Fibrational Perspectives on Determinization of Finite-State Automata}
\def\authorrunning{Thea Li}

\begin{document}
\maketitle


\begin{abstract}
  Colcombet and Petri\c{s}an argued that automata may be usefully considered from a
  functorial perspective, introducing a general notion of "$\mathcal{V}$-automaton"
  based on functors into $\mathcal{V}$.
  This enables them to recover different standard notions of automata by
  choosing $\mathcal{V}$ appropriately, and they further analyzed the
  determinization for \Rel-automata using the Kleisli adjunction between $\Set$
  and \Rel.
  In this paper, we revisit Colcombet and Petri\c{s}an's analysis from a fibrational
  perspective, building on Melliès and Zeilberger's recent alternative but
  related definition of categorical automata as functors $p : \mathcal{Q} \to \mathcal{C}$
  satisfying the finitary fiber and unique lifting of factorizations property.
  In doing so, we improve the understanding of determinization in three regards:
  Firstly, we carefully describe the universal property of determinization in
  terms of forward-backward simulations.
  Secondly, we generalize the determinization procedure for \Rel automata using a
  local adjunction between \SpanSet and \Rel, which provides us with a
  canonical forward simulation.
  Finally we also propose an alternative determinization based on the multiset
  relative adjunction which retains paths, and we leverage this to provide
  a canonical forward-backward simulation.

\end{abstract}

\section{Introduction}

%

One motivation for studying automata from a categorical point of view
is that it enables us to reason about properties of automata abstractly.
In particular, this is useful when considering constructions where canonicity
might be desirable, such as for determinization of non-deterministic automata.
Commonly one considers automata as algebras or coalgebras, see among
others \cite{adamek,EILENBERG,RUTTEN,ARBIB}.

A different approach is taken by Colcombet and Petri\c{s}an
\cite{colcombet}. They introduced a general notion of $\mathcal{V}$-automaton based
on functors into $\mathcal{V}$, which enabled them to recover standard notions of
automata by instantiating $\mathcal{V}$ appropriately.
In particular, they show that non-deterministic automata can be seen as \Rel-automata,
and they further provide a determinization for \Rel-automata using
the Kleisli-adjunction between $\Set$ and \Rel.
The \Rel-automata approach to non-deterministic automata is practical for
studying logical aspects of automata, yet it does not retain some of the
computational content.
Consequently, we revisit Colcombet and Petri\c{s}an analysis of determinization
from the point of view of Mellies and Zeilberger's \cite{mellies} alternative approach,
using functors with the \emph{unique lifting of factorizations} (ULF) property,
that ameliorates the loss of computational content and generalizes \Rel-automata.
Their approach describes non-deterministic automata as ULF-functors over a
category of labels, with finite fibers.

Both \cite{mellies} and \cite{colcombet} provide a definition of deterministic
automata, which turn out to be equivalent under the Grothendieck
construction.
As mentioned, determinization categorically has previously been studied by
Colcombet and Petri\c{s}an for \Rel-automata, but is yet to be done Melliès and
Zeilberger's approach.
In this work, we detail the determinization procedure for \Rel and
provide a generalization of this to the setting of
\cite{mellies}, recovering a similar universal property.
We then discuss the links
between our construction and the usual algorithmic procedure first presented by
Rabin and Scott \cite{scott}.
As this procedure identifies paths with the same label, which is not always
desirable, we construct an alternative determinization using multisets that
is path-relevant.

\subsection{Outline}
We recall in \cref{sec:preli} the two definitions of non-deterministic automata
as functors
and the common definition of a deterministic automaton.
In \cref{sec:reldeterminization}
we revisit the determinization procedure for \textbf{Rel}-automata and in
\cref{sec:Determinization} we construct
a generalization of this procedure to the \textbf{Span(Set)} case.
We then construct an alternative determinization in \cref{alternative} with a
stronger universal property, though the downside of this one is a significant
increase in size.
\section{Preliminaries}
\label{sec:preli}

\subsection{Automata as functors}
In \cite{colcombet}, Colcombet and Petri\c{s}an consider functors as generalized
automata. They view the domain as the category of input words and the codomain
as the category of outputs, and the functor itself as specifying the behavior of
the automaton. Additionally, to specify the accepted input, they require a full
subcategory of the category of inputs, that in practice amounts to choosing
initial and final states.

\begin{definition}\cite[Def. 2.1]{colcombet}
  \label{funct}
  A $\mathcal{V}$-automaton is a functor $\mathcal{A} : \mathcal{I}  \to \mathcal{V}$,
  together with a full subcategory inclusion $i : \mathcal{O} \to \mathcal{I}$
\end{definition}

This definition has the advantage to be very simple categorically and easy
to manipulate. Moreover, the accepted language can be defined naturally
as the composition of the automata and its full subcategory inclusion.

\begin{definition}\cite[Def. 2.1]{colcombet}
  The language accepted by a $\mathcal{V}$-automaton
  $\mathcal{A}: \mathcal{I} \to \mathcal{V}$ is the composition
  $\mathcal{A} \circ i : \mathcal{O} \to \mathcal{V}$.
\end{definition}


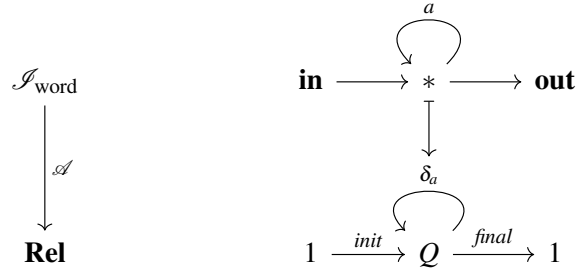
\begin{figure}
  \centering
  \begin{tikzcd}
    \mathcal{I}_{\word}\ar[dd, "\mathcal{A}"] & \quad & \mathrm{\textbf{in}}\arrow[r] &
    * \arrow[loop,looseness=6,swap,
    "a"]\arrow[r,]\arrow[d,mapsto] & \mathrm{\textbf{out}} \\
    & & & {} &\\
    \Rel & \quad & 1\arrow[r,"\emph{init}"] &
    Q\arrow[loop,looseness=4, swap,
    "\delta_a"]\arrow[r,"\emph{final}"] & 1\\
  \end{tikzcd}
  \caption{a non-deterministic automata as a functor into \Rel}
  \label{fig:Relaut}
\end{figure}

To model classical non-deterministic automata without epsilon transitions they
particularize the definition to the case where $\mathcal{V}$ is \Rel,
the category of sets and relations.
Given a fixed automaton over the alphabet $\Sigma$, they take $\mathcal{I}$ to be
the one object category associated to the free monoid $\Sigma^*$, to which they add
two objects "\textbf{in}" and "\textbf{out}" and two morphisms that act as begin-
and end-input markers. $\mathcal{O}$ is then the full subcategory generated by in and out,
we denote this by $\mathcal{I}_{\word}$.
The functor maps the base object of the free monoid to the set of states, each
letter to the relation generated by the transition function, "\textbf{in}" and
"\textbf{out}" are both mapped to 1 and the begin- and end-input markers to
relations pointing out initial and final states, see \cref{fig:Relaut}.
%

In the remainder of the paper, we will call automata with codomain \Rel,
\Rel-automata, regardless of domain.

\subsection{Non-deterministic automata as ULF functors}

The \Rel-automata approach to non-deterministic automata is
language oriented and models very well the logical aspects of automata.
Yet, as the functor is mapping into \Rel, it only interprets transitions
as ways of relating states under a label and does not account for which states
it passed along the way, nor does it account for the different ways of doing so.
It thus equates different paths, and we loose some of the computational
content of automata.

Mellies and Zeilberger \cite{mellies} introduced a categorification of
non-deterministic automata (without epsilon transitions) using ULF-functors,
that is closer to the classical definition, ameliorates the loss of computational
content, and is actually a generalization of \Rel-automata.
The idea is for the domain to be the underlying directed graph of the automaton
and the codomain to be the alphabet, with the ULF functor specifying the label
of a path.
The ULF-property then describes an important feature of automata, that a path
$\alpha$ labelled by a word $w$ can be factored in subpaths along the letters of
$w$.

\begin{definition}[ULF-functor]
  A functor $p$ is ULF if for each factorization of a morphism $p(\alpha) = uv$
  there are unique morphisms $\beta$ and $\gamma$ such that $\alpha=\beta\gamma$ and
  $p(\beta) = u$ and $p(\gamma) = v$.
\end{definition}

\begin{definition}[ULF-automaton {{\cite[Def. 3.3]{mellies}}}]
  \label{def:ulfautomaton}
  An automaton over a category $\mathcal{C}$ is a tuple $(\mathcal{C},
  \mathcal{Q} , p: \mathcal{Q} \to \mathcal{C}, q_0, Q_f)$ where $p$ is a ULF functor with
  finite fibers, $q_0$ is an object of $\mathcal{Q}$, and $Q_f$ is a set of
  objects in $\mathcal{Q}$.
\end{definition}
\noindent The objects of $\mathcal{Q}$ correspond to states, $q_0$ to the
initial state, and $Q_f$ to the final states, while unfactorizable morphisms in
$\mathcal{C}$ correspond to letters of our alphabet, and consequently arbitrary
morphism to words. This naturally leads us to consider the following as its
language.

\begin{definition}\cite[Def. 3.3]{mellies}
  The language recognized by a ULF-automaton
  $M$ with labelling functor $p$, is
  $L_M := \{ w\ |\ \exists q_f \in Q_f,\ \exists \alpha : q_0 \to q_f,\ p(\alpha)=w \}$.
\end{definition}



\begin{figure}
  \begin{align*}
    \scalebox{0.7}{
    \begin{tikzpicture}[->,shorten >=1pt,auto,node distance=2.5cm, semithick,initial text={}]
        \node[initial, state] (A)  {$1$};
        \node[state, accepting] (B) [right of=A] {$2$};
        \path (A) edge [bend left]        node {$a$} (B);
        \path (A) edge [loop above]       node {$a$} (A);
        \path (A) edge [swap, bend right] node {$b$} (B);
        \path (B) edge [loop below]       node {$b$} (B);
    \end{tikzpicture}}
  &&
  \scalebox{0.7}{
    \begin{tikzpicture}[->,shorten >=1pt,auto,node distance=2cm, semithick,initial text={}]
      \node[state] (A)  {1};
      \node[state, accepting] (B) [right of=A] {2};
      \path (A) edge [bend left, color=RoyalBlue]        node { } (B);
      \path (A) edge [loop above, color=RoyalBlue]       node { } (A);
      \path (A) edge [swap, bend right, color=Maroon] node { } (B);
      \path (B) edge [loop below, color=Maroon]       node { } (B);
      \node[state] (C) [below= 4cm of $(A)!0.5!(B)$]{};
      \path (C) edge [loop above, color=RoyalBlue]     node {$a$} (C);
      \path (C) edge [loop below, color=Maroon]       node {$b$} (C);
      \draw [shorten >=1.2cm,shorten <=0.8cm, ->] ($(A)!0.5!(B)$) -- node[left] {$p$} (C);
  \end{tikzpicture}
  }
  \end{align*}
  \caption{A non-deterministic automaton and its ULF-automaton interpretation
  (where color of the arrow indicates functorial assignment)}
  \label{fig:example}
\end{figure}
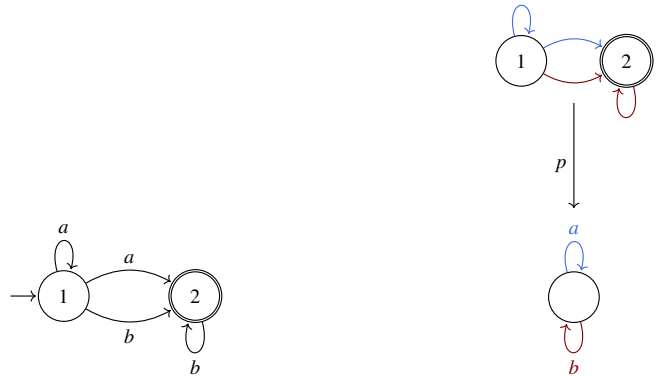

It turns out that ULF-automata generalizes \Rel-automata, indeed any
functor $\mathcal{I} \to \Rel$ can be extended to a functor
$\mathcal{I} \to \SpanSet$. Yet by a Grothendieck-like correspondence
any pseudofunctor into \SpanSet corresponds to a ULF-functor
\cite{mellies}, hence any \Rel-automaton induces a ULF-automaton
over $\mathcal{I}$. We can in fact characterize these as a particular type of
ULF-functors.

\begin{proposition}
    A functor $F : \mathcal{C} \to \SpanSet$ factors through
    \Rel if and only if its corresponding functor lying over
    $\mathcal{C}$ (by the Grothendieck construction) is ULF and
    faithful.
\end{proposition}

\begin{corollary}
  \Rel-automata $\mathcal{C} \to \Rel$ factoring through \textbf{FinRel} are
  equivalent to faithful ULF-automata over $\mathcal{C}$.
\end{corollary}

We can, as noted in \cite{mellies}, represent ULF-automata as pseudofunctors
$\mathcal{C} \to \SpanSet$ factoring through \SpanFin, up to choice of initial and final
states.
In practice, we will often directly work with the equivalent
\SpanSet-functor rather than with the ULF-automaton, as it makes it
easier to relate automata over the same category.
We will henceforth refer to them as \SpanSet-automata.

We remark that we on occasion drop the condition that our
\SpanSet-automata factor through \SpanFin as it is not always important.
In \cref{sec:Determinization} all constructions preserve finiteness, and
consequently the construction naturally restricts to \SpanFin, thus we
will consider general \SpanSet-automata.
However, in \cref{alternative} whether constructions are finite is important
and we adapt our definitions accordingly.


\subsection{Deterministic Automata}
Both \cite{colcombet} and \cite{mellies} give definitions of a deterministic
automaton. These definitions turn out to be equivalent under the Grothendieck
construction given the right instantiations, and being able to reason using both
approaches is useful.

By definition an automaton is deterministic when its transition relation is a
functional relation, consequently, \cite{colcombet} models deterministic
automata as \Rel-automata that factor through $\Set$.

\begin{definition}[\cite{colcombet}]
  A deterministic automaton is a functor $\mathcal{I}_{\word} \to \Set$,
  which maps \emph{in} to 1 and \emph{out} to 2.
\end{definition}

For a ULF-functor to model a deterministic automaton it not only needs a unique
lift of factorizations of words, as noted in \cite{zeilberger}, it requires
every word to have a unique path given a starting state, meaning that any word
$w: p(q) \to c$ has a unique lift $\alpha : q \to q'$.

\begin{definition}[Discrete opfibration]
  A functor $p : \mathcal{Q} \to \mathcal{C}$ is a discrete opfibration if any
  morphism $w : p(q) \to c$ has a unique lift, i.e. there exists a unique
  $\alpha : q \to q'$ such that $p(\alpha) = w$.
\end{definition}

\begin{definition}\cite[Def. 2.9]{zeilberger}
  A deterministic automaton over a category $\mathcal{C}$ is a tuple
  $(\mathcal{C}, \mathcal{Q}, p: \mathcal{Q} \to \mathcal{C}, q_0, Q_f)$ where
  $p$ is a discrete opfibration with finite fibers, $q_0$ is an object of
  $\mathcal{Q}$, and $Q_f$ is a set of objects in $\mathcal{Q}$.
\end{definition}

It turns out that these two definitions are equivalent given the right instantiations.

\begin{proposition}
  By the Grothendieck construction, functors $\mathcal{C} \to \Set$ are
  equivalent to discrete opfibrations over $\mathcal{C}$.
\end{proposition}

\begin{corollary}
  $\Set$-automata $\mathcal{I}_{\word} \to \Set$ are equivalent to
  deterministic automata over $\mathcal{I}_{\word}$ where the fiber over
  $\mathrm{in}$ has one element and the fiber over $\mathrm{out}$ has two.
\end{corollary}

\section{Determinisation for Rel}
\label{sec:reldeterminization}
In order to determinize \SpanSet-automata, we first focus on the particular case
of \Rel-automata, previously studied in \cite{colcombet}, this procedure will
constitute in one building stone for our \SpanSet-automata determinization.
We recall the determinization procedure for \Rel-automata in \cref{subsec:detrel},
and transpose it to faithful ULF-automata.
Finally we discuss the difficulties of extending the construction to
\SpanSet.



\subsection{Determinisation for Rel}
\label{subsec:detrel}

Determinizing a \Rel-automaton $\mathcal{I} \to \Rel$ consists in
transforming it into a functor $\mathcal{I} \to \Set$; to obtain such,
one can postcompose by a functor $\Rel \to \Set$.
As explained in \cite{colcombet}, such a functor naturally arises, as
\Rel is the Kleisli category of the powerset monad on $\Set$.
Consequently, there is an adjunction
${L : \Set \rightleftarrows \Rel : R}$.
The functor $R$ sends objects $A$ in \Rel to their powersets
$\mathcal{P}(A)$ and relations $r \subseteq A \times B$ to the functions
$\mathcal{P}(A)\to \mathcal{P}(B)$ that sends a subset $S\subseteq A$ to the
subset of all elements in $B$ related to some element in $S$.

\begin{construction}\cite[Sec 3.3]{colcombet}
    Given a \Rel-automaton $F:\mathcal{I} \to \Rel$, its
    determinization is given by
    \begin{equation*}
        Det(F) :=
        \begin{tikzcd}
            \mathcal{I} \arrow[r, "F"] & \Rel \arrow[r, "R"] & \Set
        \end{tikzcd}
    \end{equation*}
\end{construction}

To be an appropriate determinization procedure the deterministic automaton
should simulate the original automaton and recognize the same language.
Simulations can be represented categorically as the following, which is inspired
by simulations of labelled transition systems presented in \cite{sobociski} and
a slight modification to the morphisms of automata presented in \cite{colcombet}.
Moreover, this gives a suitable notion of morphism of automata.

\begin{definition}
  Given two automata $F, F': \mathcal{I} \to \Rel$,
  a (forward-backward) \emph{simulation} from $F$ to $F'$ is a natural
  transformation $\alpha : F' \Rightarrow F$.
\end{definition}
Indeed for any $c, d \in \text{ob}(\mathcal{I})$ and any word $w:c\to d$ if
$(q', s)\in F(w)\circ \alpha_c$ then there is a state $q \in F'(c)$ such that
$(q,q')\in F(w)$ and $(q,s)\in \alpha_c$, by naturality we also have that
$(q', s)\in \alpha_d \circ F'(w)$ consequently there is a state $s'\in F(d)$
and $(s,s')\in F'(w)$ and $(q',s')\in \alpha_d$, that is, $F'$ simulates $F$.
We can make the reversed argument for any $(q',s)\in \alpha_d \circ F'(w)$, that
is, we have a backward simulation.

\begin{definition}
    \Rel-automata with domain $\mathcal{I}$ and simulations
    form a category Aut($\mathcal{I}$).
\end{definition}


We can then justify the canonicity of the determinization by the existence of a
canonical (\emph{forward-backwards}) simulation from an automaton to its
determinization, in which the initial state is simulated by an initial state and
the final states are simulated by final states.

\begin{proposition}
    Let $F: \mathcal{I}\to \Rel$ be an automaton,
    then there is a canonical (forward-backward) simulation
    from $F$ to Det($F$) given by the counit of the Kleisli-adjunction
    between \Rel and $\Set$.
    \begin{equation*}
        \begin{tikzcd}
            \mathcal{I} \arrow[rr, "F"] &  & \Rel \ar[equal]{rr} \arrow[rd, "R"'] & {} & \Rel \\
                                          &  && \Set \arrow[ru, "i"'] \arrow[u, "\varepsilon"', Rightarrow] &
        \end{tikzcd}
    \end{equation*}
\end{proposition}

Having a canonical simulation is not enough to justify that the
determinization preserves the language, as the simulation is one-sided, and as
we have yet to address how initial and final states are determined for Det($F$).
Under Grothendieck construction Det($F$) is
\[ \text{Det}(M):=(\int R F, \mathcal{I}, \Pi:\int R F\to \mathcal{I},
q'_0, Q'_f) \]

We take as initial state $q'_0 = (p^{-1}(c), \{q_0\})$, where $p$ is the
ULF-functor corresponding to $F$, $q_0$ is the initial state of $F$ and $c =
p(q_0)$.
The set of final states $Q'_f$ are all states $(p^{-1}(q_f), S)$ such that $q_f$
is a final state in $F$ and $q_f\in S$.
Given this it is straightforward to check that Det($F$) recognizes the same
language as $F$.

\subsection{Universal Property}

Now that we have justified the canonicity of our determinization procedure, we
would like to be able to relate it to other determinizations and other
deterministic automata, that is, equip it with a universal property.
To characterize the universal property we make use of the notion of bisimulation.
For this purpose, recall that \Rel is a dagger category taking
$(-)^\dagger$ as the converse relation.


\begin{definition}
    If $\alpha : F' \Rightarrow F$ is a natural transformation
    such that $(\alpha)^\dagger$ also is a natural transformation, then
    $\alpha$ is called a \emph{bisimulation}.
\end{definition}

Indeed, requiring $(\alpha)^\dagger$ to also be a simulation gives us the converse
simulation, i.e. "$F$ simulates $F'$", by the simulation $\alpha$.

As our canonical simulation is built out of a counit we make use of the triangle
identity for unit and counit of adjoint functors to obtain the following
universal property.

\begin{proposition}\label{UniRel}
    Let $G:\mathcal{I} \to \Set$ be some deterministic automaton then any
    \emph{forward-backward} simulation from $F$ to $G$ factors uniquely through the
    universal simulation to Det($F$) and a bisimulation between Det($F$) and
    $G$.
\end{proposition}
\begin{proof}
    By one of the triangle identities for adjoint functors we have
    the following equality.
    \begin{equation*}
        \begin{tikzcd}
            \mathcal{I} \arrow[rr, "F"{name=Fp}] \arrow[rd, "G"']
                &
                & \Rel \\
                & \Set \arrow[ru, "i"']
                               \arrow[u, to=Fp, shorten >=0.7ex, "\alpha"', Rightarrow]
                &
        \end{tikzcd}
        =
        \begin{tikzcd}
            \mathcal{I} \arrow[rr, "F"] \arrow[rd, "G"']
                & {}
                & \Rel \arrow[equal, ""{name=eq2, below}]{rr}
                               \arrow[rd, "R"'] & {}
                & \Rel \\
                & \Set \arrow[ru, "i"']
                               \arrow[equal, ""{name=eq}]{rr}
                               \arrow[u, swap, "\alpha", Rightarrow]
                & \ar[u, swap, Rightarrow, from=eq, "\eta"]
                & \Set \arrow[ru, "i"'] \arrow[u, to=eq2, "\varepsilon"', Rightarrow]
                &
        \end{tikzcd}
    \end{equation*}
    The simulation between Det($F$) and $G$ is a bisimulation as
    any natural transformation generated by a natural transformation
    between functors into $\Set$ and an inclusion into \Rel.
\end{proof}
%

\subsection{Difficulties in generalization to \SpanSet-automata}

This determinization procedure has the advantage to be direct and fairly simple,
but does not scale easily to \SpanSet as it is less well behaved than \Rel.
First, we do not have an adjunction between \SpanSet and $\Set$, which was
essential in the \Rel-case for determinization and to establish the universal
property.
Second, as morphisms in \SpanSet are composed using pullbacks,
functors out of or to \SpanSet tend to be lax or pseudo as pullbacks
are only unique up to isomorphism.


\section{Powerset determinization for \SpanSet-automata}
\label{sec:Determinization}

In this section, we construct local adjunction between \SpanSet and \Rel seen as
2-categories, and we show that this is sufficient to recover a determinization
procedure and universal property analogous to the \Rel-determinization for
ULF-automata/\SpanSet-automata.
We then discuss the link between the categorical determinization procedure and
the algorithmic one.

\subsection{Determinization procedure}

As ULF-automata over $\mathcal{C}$ are equivalent to pseudofunctors $\mathcal{C}
\to \SpanSet$ and deterministic ULF-automata are equivalent to functors
$\mathcal{C}\to \Set$, it suffices to be able to build a functor $\mathcal{C}\to
\Set$ out of a pseudo-functor $\mathcal{C} \to \SpanSet$ to determinize
ULF-automata.
As we know how to determinize \Rel-automata, it is sufficient for us to be
able to turn a pseudo-functor $\mathcal{C} \to \SpanSet$ into a functor
$\mathcal{C}\to \Rel$ in a canonical way.
That is, it suffices to construct a lax functor $\SpanSet \to \Rel$ that strictifies
the pseudo-functor to a strict functor $\mathcal{C}\to \Rel$ by postcomposition.
To do so, we must work in a 2-categorical setting as we deal with pseudo-functors,
and view \SpanSet and \Rel as a 2-categorical construction.

\begin{definition}
\SpanSet is a 2-category with 2-morphisms given by morphisms of spans

$A \xleftarrow{f} S \xrightarrow{g} B$ to $A \xleftarrow{f'} S' \xrightarrow{g'} B$
which are functions $h : S \to S'$ such that:
\begin{center}
    \begin{tikzcd}
        & S \arrow[d, "h"] \arrow[ld, "f"'] \arrow[rd, "g"] &   \\
      A & S' \arrow[l, "f'"] \arrow[r, "g'"']               & B
    \end{tikzcd}
\end{center}
\end{definition}
In this context, we can view \Rel as a 2-category in a similar way,
2-morphisms would be given by inclusions of relations, which indeed is just a
morphism of spans.

There is a "forgetful" lax functor from \SpanSet to \Rel
that is identity on objects and maps spans to the relation
generated by the image of the legs, we call this functor Im for
image.
As this functor identifies spans that generate the same relation,
by post-composing it to a pseudo-functor such as our automaton $F$ we obtain a
functor into \Rel that is strict.
Post-composing an automaton $F: \mathcal{C}\to
\SpanSet$ with $R \circ \text{Im}$ does indeed
give us a determinization of $F$ that preserves the language,
as the image functor only identifies all paths with the same label
between two given states.

\begin{construction}
    Given a \SpanSet-automaton $F:\mathcal{C}\to \SpanSet$, its
    determinization is given by
    \begin{equation*}
        Det(F) :=
        \begin{tikzcd}
            \mathcal{C} \arrow[r, "F"] & \SpanSet \arrow[r, "\Im"] & \Rel \arrow[r, "R"] & \Set
        \end{tikzcd}
    \end{equation*}
\end{construction}


As for \Rel-automata, for our determinization to be appropriate, the
deterministic automaton must simulate the original automaton and recognize the
same language. This requires us to adapt our notion of simulation
to the 2-categorical setting and to \SpanSet.

\begin{definition}\label{def:sim}
    Given two automata $F, F': \mathcal{C}\to \SpanSet$,
    a \emph{(forward) simulation} from $F$ to $F'$ is a lax natural
    transformation $\alpha : F' \Rightarrow F$.
\end{definition}

In relaxing our notion of simulation to a \emph{lax} natural transformation we have a
simulation by a similar reasoning as for \Rel-automata by invoking the laxness,
but not a backward simulation, as the lax naturality only guarantee a
span-morphism in one direction. If such a simulation is pseudo we do obtain the
backward simulation as well.

\begin{definition}\label{def:sim2}
  Given two automata $F, F': \mathcal{C}\to \SpanSet$,
  a \emph{(forward-backward) simulation} from $F$ to $F'$ is a pseudo natural
  transformation $\alpha : F' \Rightarrow F$.
\end{definition}

\begin{definition}
    \SpanSet-automata with domain $\mathcal{C}$ and forward
    simulations form a category Aut($\mathcal{C}$).
\end{definition}

Though the image functor does not have any strict adjoints, it is the left
adjoint of a local adjunction between \Rel and \SpanSet seen
as 2-categories, with the right adjoint being the inclusion.
Using this local adjunction properly will enable us the prove the canonicity of our
determinization procedure, as well as recover a universal property similar to
one of \Rel. We start by recalling the definition and terminology used for local
adjunctions.

\begin{definition}
  Let $\mathcal{C}$ and $\mathcal{D}$ be bicategories, a local adjunction between $\mathcal{C}$
  and $\mathcal{D}$ is comprised of two functors $L:\mathcal{C}\to\mathcal{D}$ and
  $R:\mathcal{D}\to\mathcal{C}$, inducing a family of adjunctions
  \begin{center}
    \adjustbox{scale=0.7}{
      $\begin{tikzcd}
      \mathcal{D}(LA,B) \arrow[rr, bend left] & \bot & \mathcal{C}(A,RB) \arrow[ll, bend left]
      \end{tikzcd}$}
  \end{center}
  natural in $A$ and $B$. We use $\eta$ and $\epsilon$ to denote units and counits
  on the local level.
\end{definition}

To construct a canonical simulation we begin by constructing a canonical
simulation from a \SpanSet-automaton and its \Rel-ification,
for this purpose we will make use of the following definition.

\begin{definition}
    We denote by $"\eta" : i\Im \Rightarrow \Id_{\SpanSet}$ the lax natural
    transformation that is the identity on each component with lax naturality
    given by the $\eta$ the unit of the local adjunction
    $\Im : \SpanSet(A,B) \rightleftarrows \Rel(A,B) : i$.
    We respectively denote by $"\epsilon": \Id_{\Rel} \Rightarrow \Im i$ the
    lax natural transformation associated to the counit. Let $c$ and $d$ be sets
    and $w$ and and $R$ be a span and a relation from $c$ to $d$, respectively,
    then the lax naturality squares are the following for $"\eta"$ and $"\epsilon"$:
    %
    \begin{align*}
        \begin{tikzcd}[ampersand replacement = \&]
            c \arrow[rr, "i(\text{Im}(w))"] \arrow[d, "\Id_c"']
            \& {} \& d \arrow[d, "\Id_d"] \\
            c \arrow[rr, "w"']   \arrow[rru, Rightarrow, "{\eta}", shorten =15pt]
            \& {} \& d
        \end{tikzcd}
        &&
        \begin{tikzcd}[ampersand replacement = \&]
            c \arrow[rr, "R"] \arrow[d, "\Id_c"']
            \& {} \& d \arrow[d, "\Id_d"] \\
            c \arrow[rr, "\Im (i (R))"']   \arrow[rru, Rightarrow, "{\epsilon}", shorten =15pt]
            \& {} \& d
        \end{tikzcd}
    \end{align*}
\end{definition}

\begin{lemma}\label{SpanRel}
    Let $F: \mathcal{C}\to \SpanSet$ be an automaton,
    then there is a canonical (forward) simulation from $F$ to
    $i \text{Im} F$ given by the lax natural transformation $"\eta"$ defined
    above.
    \begin{equation*}
        \begin{tikzcd}[ampersand replacement = \&]
            \mathcal{C} \arrow[rr, "F"] \& \& \SpanSet \arrow[rr, equal] \arrow[rd, "\Im"']
            \& {}
            \& \SpanSet \\
            \& {} \& {} \& \Rel \arrow[u, "{"\eta"}"', Rightarrow] \arrow[ru, "i"']
            \&
        \end{tikzcd}
    \end{equation*}
\end{lemma}

Using this we can now obtain a canonical simulation from an automaton to its
determinization similar to the one for \Rel.

\begin{proposition}
    Let $F: \mathcal{C}\to \SpanSet$ be an automaton,
    then there is a canonical (forward) simulation
    from $F$ to Det($F$) given by pasting the (forward) simulation
    \ref{SpanRel} and the counit of the Kleisli-adjunction between \Rel and $\Set$.
    \begin{equation*}
        \adjustbox{scale=0.8}{
        \begin{tikzcd}
            \mathcal{C} \arrow[rrr, "F"]
            & {} &
            & \SpanSet \arrow[rr, equal] \arrow[rd, "\Im"]
            & {} & \SpanSet \arrow[rr, equal] & & \SpanSet \\
            & &
            & {}
            & \Rel \arrow[ru, "i"'] \arrow[rr, equal] \arrow[rd, "R"] \arrow[u, "{"\eta"}"', Rightarrow]
            & {} & \Rel \arrow[ru, "i"'] & \\
            &
            & {} &
            & & \Set \arrow[ru, "i"'] \arrow[u, "\varepsilon"', Rightarrow]
            & &
            \end{tikzcd}}
    \end{equation*}
\end{proposition}

Moreover, it indeed recognizes the same language as postcomposing with \text{Im}
does not affect the language and the analogous result for \Rel-automata.
Hence this determinization procedure is a fitting determinization procedure.

\subsection{Universal Property for \SpanSet}

Having proven the canonicity of our determinization, we now construct the
universal property of the determinization of \SpanSet-automata.
As before, we rely on bisimulations, which we need to adapt to \SpanSet.
For this purpose, recall \SpanSet is a dagger category.

\begin{proposition}
    \SpanSet is a dagger category taking $(-)^\dagger$ as the converse
    span.
\end{proposition}
\begin{definition}
    A a natural transformation $\alpha : F' \Rightarrow F$ such that
    $(\alpha)^\dagger : F \Rightarrow F'$ is also a natural transformation is a
    \emph{bisimulation}.
\end{definition}

Note that any such alpha must incidentally be strict or pseudo.

To recover a similar universal property as for \Rel, we additionally need a
pasting inverse to the simulation $"\eta"$.
To construct such an inverse we leverage the fact that our construction is based
on a local adjunction.

\begin{lemma}\label{lemma:inverseSpan}
    The (forward) simulation $"\epsilon"$ is the left pasting inverse of the forward
    simulation
    $"\eta"$
    \[
\begin{array}{ccc}
    \begin{tikzcd}[ampersand replacement=\&]
        \& \SpanSet \arrow[rr, equal] \arrow[rd, "\text{Im}"']  \& {} \& \SpanSet \\
        \Rel \arrow[ru, "i"] \arrow[rr, equal] \& {} \arrow[u, shorten <=0.8ex,  "{"\varepsilon"}" description, Rightarrow] \& \Rel \arrow[u, "{"\eta"}" description, Rightarrow] \arrow[ru, "i"']   \&
    \end{tikzcd}
    &
    =
    &
    \begin{tikzcd}[ampersand replacement=\&]
        \SpanSet           \\
        \Rel \arrow[u, bend left=49, "i"{name=i}] \arrow[u, bend right=49, "i"'{name=i2}] \arrow[from=i, to=i2, Rightarrow, shorten =7pt, "\textbf{id}_i"']
    \end{tikzcd}
\end{array}\]
\end{lemma}
Note that as $"\epsilon"$ is the identity, this lemma says that
$"\eta"$ restricted to \Rel is the identity natural transformation.

Using this "triangle-identity" for $"\epsilon"$ and $"\eta"$ and the universal
property for \Rel-automata we can now obtain a universal property for
\SpanSet-automata.

\begin{proposition}
    Let $G:\mathcal{C}\to\Set$ be some deterministic
    automaton, then any forward simulation from $F$
    to $G$ factors through the canonical simulation to Det($F$)
    and a bisimulation between Det($F$) and $G$.
\end{proposition}

\begin{proof}
    By \cref{UniRel} and \cref{lemma:inverseSpan} we have that
    the following equality of diagrams:
    \begin{equation*}
        \adjustbox{scale=0.65}{
        \begin{tikzcd}
            \mathcal{C} \arrow[rdd, "G"', bend right] \arrow[rrr, "F"]
            & {} &
            & \SpanSet\\
            & & \Rel \arrow[ru, "i"]\\
            & \Set \arrow[ru, "i"]\arrow[uu, "\alpha"', Rightarrow]
            \end{tikzcd}}
            =
            \adjustbox{scale=0.65}{
                \begin{tikzcd}
                    \mathcal{C} \arrow[rdd, "G"', bend right] \arrow[rrr, "F"]
                    & {} &
                    & \SpanSet \arrow[rr, equal] \arrow[rd, "\text{Im}"]
                    & {} & \SpanSet \arrow[rr, equal] & & \SpanSet \\
                    & & \Rel \arrow[ru, "i"] \arrow[rd, "R"] \arrow[rr, equal]
                    & {} \arrow[u, shorten <=0.8ex, "{"\varepsilon"}"', Rightarrow]
                    & \Rel \arrow[ru, "i"'] \arrow[rr, equal] \arrow[rd, "R"] \arrow[u, "{"\eta"}"', Rightarrow]
                    & {} & \Rel \arrow[ru, "i"'] & \\
                    & \Set \arrow[ru, "i"] \arrow[rr, equal] \arrow[uu, "\alpha"', Rightarrow]
                    & {} \arrow[u, shorten <=0.8ex, "\eta"', Rightarrow] & \Set \arrow[rr, equal]
                    & & \Set \arrow[ru, "i"'] \arrow[u, "\varepsilon"', Rightarrow]
                    & &
                \end{tikzcd}}
    \end{equation*}
    Consequently this gives a unique factorization of $\alpha$ through
    the universal simulation between $F$ and Det($F$).
\end{proof}

\subsection{Relation to the classical determinization algorithm}

To see how this version of determinization presented here
relates to the classical determinization algorithm we start
by working through an example.

\begin{example}
    Consider the following automaton with its ULF-automaton interpretation:
    \begin{align*}
    \scalebox{0.7}
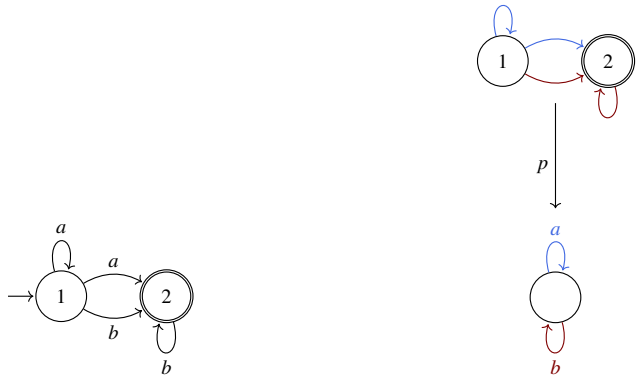
{
        \begin{tikzpicture}[->,shorten >=1pt,auto,node distance=2cm, semithick,initial text={}]
            \node[initial, state] (A)  {$1$};
            \node[state, accepting] (B) [right of=A] {$2$};
            \path (A) edge [bend left]        node {$a$} (B);
            \path (A) edge [loop above]       node {$a$} (A);
            \path (A) edge [swap, bend right] node {$b$} (B);
            \path (B) edge [loop below]       node {$b$} (B);
        \end{tikzpicture}}
    &&
    \scalebox{0.7}{
    \begin{tikzpicture}[->,shorten >=1pt,auto,node distance=2cm, semithick,initial text={}]
      \node[state] (A)  {1};
      \node[state,accepting] (B) [right of=A] {2};
      \path (A) edge [bend left, color=RoyalBlue]        node { } (B);
      \path (A) edge [loop above, color=RoyalBlue]       node { } (A);
      \path (A) edge [swap, bend right, color=Maroon] node { } (B);
      \path (B) edge [loop below, color=Maroon]       node { } (B);
      \node[state] (C) [below= 4cm of $(A)!0.5!(B)$]{};
      \path (C) edge [loop above, color=RoyalBlue]     node {$a$} (C);
      \path (C) edge [loop below, color=Maroon]       node {$b$} (C);
      \draw [shorten >=1.2cm,shorten <=0.8cm, ->] ($(A)!0.5!(B)$) -- node[left] {$p$} (C);
  \end{tikzpicture}}
    \end{align*}
    Using the determinization process presented in the previous section, the
    determinization of this automaton is the ULF-automaton pictured below to the
    right, which is equivalent to what we would get using the algorithmic
    procedure, pictured below to the left.

    \begin{align*}
        \scalebox{0.7}
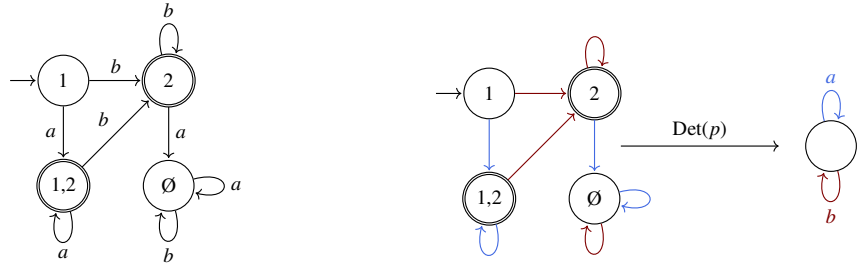
{
            \begin{tikzpicture}[->,shorten >=1pt,auto,node distance=2cm, semithick,initial text={}]
                \node[initial, state] (A)  {1};
                \node[state, accepting] (B) [right of=A] {2};
                \node[state, accepting] (F) [below of=A] {1,2};
                \node[state] (G) [below of=B] {\O};
                \path (A) edge  [swap]    node { $a$} (F);
                \path (A) edge            node {$b$ } (B);
                \path (B) edge [loop above]       node {$b$ } (B);
                \path (F) edge [] node {$b$ } (B);
                \path (F) edge [ loop below] node {$a$ } (F);
                \path (B) edge [] node {$a$ } (G);
                \path (G) edge [loop right] node {$a$ } (G);
                \path (G) edge [loop below] node {$b$} (G);
            \end{tikzpicture}}
        &&
        \scalebox{0.7}{
        \begin{tikzpicture}[->,shorten >=1pt,auto,node distance=2cm, semithick,initial text={}]
            \node[initial, state] (A)  {1};
            \node[state, accepting] (B) [right of=A] {2};
            \node[state, accepting] (F) [below of=A] {1,2};
            \node[state] (G) [below of=B] {\O};
            \path (A) edge [color=RoyalBlue]        node { } (F);
            \path (A) edge [swap, color=Maroon] node { } (B);
            \path (B) edge [loop above, color=Maroon]       node { } (B);
            \path (F) edge [color=Maroon] node { } (B);
            \path (F) edge [color=RoyalBlue, loop below] node { } (F);
            \path (B) edge [color=RoyalBlue] node { } (G);
            \path (G) edge [color=RoyalBlue, loop right] node { } (G);
            \path (G) edge [color=Maroon, loop below] node { } (G);
            \node[state] (C) [right= 5cm of $(A)!0.5!(G)$]{ };
            \path (C) edge [loop above, color=RoyalBlue]     node {$a$} (C);
            \path (C) edge [loop below, color=Maroon]       node {$b$} (C);
            \draw [shorten >=0.5cm,shorten <=0.5cm, ->] ($(B)!0.5!(G)$) -- node[above] {Det($p$)} (C);
        \end{tikzpicture}}
        \end{align*}
\end{example}

Let us recall the algorithmic procedure for determinizing an automaton. Given an
automaton $M = (\Sigma, Q, \delta : Q \times \Sigma \to \mathcal{P}(Q), q_0, Q_f)$ we begin by taking the
powerset of the set of states $\mathcal{P}(Q)$, this will be the states of our
deterministic machine.
Then for the transition function we just take the function $\delta': \mathcal{P}(Q)\times\Sigma\to
\mathcal{P}(Q)$, generated by $\delta$, mapping $S \subseteq Q \mapsto \{q \in Q\ | \ \exists q'(q
\in \delta(q'))\}$.

When we have a free monoid as base category for our automaton the
determinization we obtain is analogous to the algorithmic
one. This is because postcomposition by $R$ replaces the objects in
$\mathcal{C}$ with the powersets of their fibers,
if we only have one fiber, then the determinization has its powerset as set of
states.

For a general base category we do not get the full powerset of the original set
of states, as the states of our deterministic machine.
Instead the determinization process takes to a certain extent into account which
subsets would be unreachable.
The limiting factor here is that these are only the subsets comprised of states
lying over different fibers and consequently have no chance of a path to or from
them.

\begin{example}
    Consider the following automaton $M$, with ULF-functor $p$:
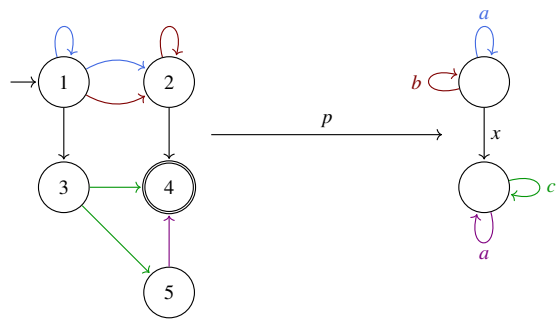
    \begin{center}
        \scalebox{0.7}{
        \begin{tikzpicture}[->,shorten >=1pt,auto,node distance=2cm, semithick,initial text={}]
            \node[initial, state] (A)  {1};
            \node[state] (B) [right of=A] {2};
            \node[state]    (X) [below of=A] {3};
            \node[state, accepting]    (Y) [below of=B] {4};
            \node[state]    (Z) [below of=Y] {5};
            \path (A) edge [bend left, color=RoyalBlue]        node { } (B);
            \path (A) edge [loop above, color=RoyalBlue]       node { } (A);
            \path (A) edge [swap, bend right, color=Maroon] node { } (B);
            \path (B) edge [loop above, color=Maroon]       node { } (B);
            \path (A) edge node { } (X);
            \path (B) edge node { } (Y);
            \path (X) edge [color=myGreen] node { } (Y);
            \path (X) edge [color=myGreen] node { } (Z);
            \path (Z) edge [color=Purple] node { } (Y);
            \node[state] (C) [right=5cm of B] { };
            \node[state] (W) [below of= C] { };
            \path (C) edge [loop above, color=RoyalBlue]     node {$a$} (C);
            \path (C) edge [loop left, color=Maroon]       node {$b$} (C);
            \path (C) edge node {$x$} (W);
            \path (W) edge [loop right, color=myGreen]     node {$c$} (W);
            \path (W) edge [loop below, color=Purple]     node {$a$} (W);
            \draw [shorten >=0.8cm,shorten <=0.8cm, ->] ($(B)!0.5!(Y)$) -- node[above] {$p$} ($(C)!0.5!(W)$);
        \end{tikzpicture}}
        \end{center}
        It recognizes the language $\{a^nb^mxcd\} \cup \{a^nb^mxc\}$,
        in particular has two "disjoint" alphabets $\{a,b\}$ and $\{c,d\}$ with
        a connective $x$, meaning that we have no paths combining them without
        an $x$ in the middle, and once we have "put" the $x$ after a possibly
        empty word $a^nb^m$ we cannot go back and only have $\{c,d\}$ available.

        The following is the determinization of $M$ using our determinization:
        \begin{center}
            \scalebox{0.7}
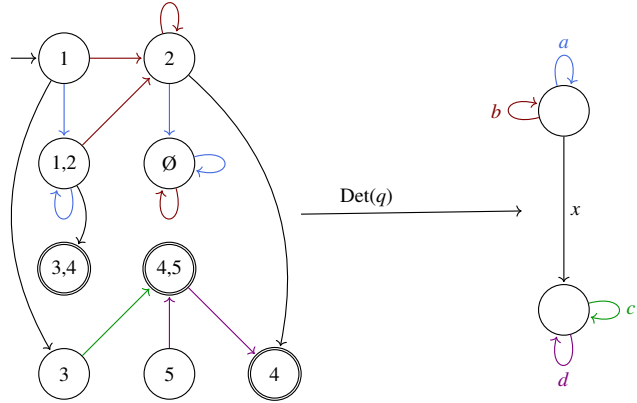
{
            \begin{tikzpicture}[->,shorten >=1pt,auto,node distance=2cm, semithick,initial text={}]
                \node[initial, state] (A)  {1};
                \node[state] (B) [right of=A] {2};
                \node[state] (F) [below of=A] {1,2};
                \node[state] (G) [below of=B] {\O};
                \path (A) edge [color=RoyalBlue]        node { } (F);
                \path (A) edge [swap, color=Maroon] node { } (B);
                \path (B) edge [loop above, color=Maroon]       node { } (B);
                \path (F) edge [color=Maroon] node { } (B);
                \path (F) edge [color=RoyalBlue, loop below] node { } (F);
                \path (B) edge [color=RoyalBlue] node { } (G);
                \path (G) edge [color=RoyalBlue, loop right] node { } (G);
                \path (G) edge [color=Maroon, loop below] node { } (G);
                \node[state] (C) [right= 8cm of $(A)!0.5!(G)$]{ };
                \path (C) edge [loop above, color=RoyalBlue]     node {$a$} (C);
                \path (C) edge [loop left, color=Maroon]       node {$b$} (C);
                \node[state, accepting] (Y) [below of=F] {3,4};
                \node[state] (X) [below of=Y] {3};
                \node[state, accepting] (Z) [below of=G] {4,5};
                \node[state] (W) [below of=Z] {5};
                \node[state, accepting] (V) [right of=W] {4};
                \path (A) edge [bend right] node { } (X);
                \path (F) edge [bend left] node { } (Y);
                \path (B) edge [bend left] node { } (V);
                \path (X) edge [color=myGreen] node { } (Z);
                \path (Z) edge [color=Purple] node { } (V);
                \path (W) edge [color=Purple] node { } (Z);
                \node[state] (U) [below= 2.8cm of C] { };
                \path (U) edge [loop right, color=myGreen]     node {$c$} (U);
                \path (U) edge [loop below, color=Purple]       node {$d$} (U);
                \path (C) edge node {$x$} (U);
                \draw [shorten >=0.8cm,shorten <=2.5cm, ->] ($(G)!0.5!(Z)$) -- node[above] {Det($q$)} ($(C)!0.5!(U)$);
            \end{tikzpicture}}
            \end{center}
            We see here that we do not get states containing states from
            different fibers, which makes our determinized automaton smaller
            that the powerset.
\end{example}

\section{Multiset determinization for \SpanSet-automata}\label{alternative}
One step in the determinization process presented
above that might feel unjustified is the identification
of paths between states having the same label (i.e. postcomposition with
the image functor into \Rel).
To make amends for this we might consider
a path-relevant version of determinization using multisets.
It is worth noting that this significantly increases the size of
our final machine, however, it strengthens the universal property.

\subsection{Procedure}
Multisets only correspond to finite spans, consequently we will restrict
ourselves to \SpanFin-automata.
This is not unjustified as we in the definition of ULF-automaton have the
condition that the fibers of the ULF-functor are finite, that is, its
corresponding \SpanSet-functor factors through \SpanFin.
To define the multiset determinization use that the multiset
functor forms a relative monad $M: \Fin \to \Set$ and the fact that
we for the "relative" Kleisli-category on \Fin, have that $\Fin_M(A,B)$ is the
skeleton of $\SpanFin(A,B)$  for finite sets $A$ and
$B$. For a thorough treatment of the theory of
relative monads we refer the reader to \cite{Altenkirch_2015,Arkor_2024}.


\begin{proposition}
  Multisets form a relative monad from \Fin to
  $\Set$, defined for all $A,B : \Fin$ and $\varphi : A \to \mathbb{N}^B$ by
  \begin{align*}
    \begin{array}{ccc}
      M : \Fin & \to     & \Set \\
      A        & \mapsto & \mathbb{N}^A
    \end{array}
    &&
    \begin{array}{ccc}
      \eta_A : A    & \to     & \mathbb{N}^A \\
      a          & \mapsto & \Sigma_{a'\in A} \chi_a{a'}
    \end{array}
    &&
    \begin{array}{ccc}
      \varphi^* : \mathbb{N}^A & \to     & \mathbb{N}^B \\
      f               & \mapsto & \Sigma_{a\in A} f(a) \cdot \varphi (a)(-)
    \end{array}
  \end{align*}
\end{proposition}

\begin{lemma}\label{localEQ}
    There is a local equivalence of categories between \SpanFin and
    $\Fin_M$ that is identity on objects defined by
    \begin{equation*}
      \begin{array}{r c l}
        U: \SpanFin(A,B)         & \rightleftarrows & \Fin_M (A,B) : \CMamalg \\
        A \xleftarrow{f} S \xrightarrow{g} B & \mapsto & (a \mapsto (b \mapsto \abs{ (f \times g)^{-1}(a,b)}))\\
        \langle \pi_1 , \pi_2 \rangle : \CMamalg_{(a,b) \in A \times B} f(a)(b) \rightrightarrows \langle A , B \rangle & \mapsfrom & f : A \to \mathbb{N}^B
      \end{array}
    \end{equation*}
\end{lemma}

\begin{construction}
  Given a \SpanFin-automaton $F:\mathcal{C}\to \SpanFin$, its multiset
  determinization is given by the following,
  where $\Fin_M$ is the relative Kleisli-category associated to $M$.
  \begin{equation*}
    \mathrm{MDet}(F) :=
    \begin{tikzcd}
      \mathcal{C} \arrow[r, "F"] & \SpanFin \arrow[r, "U"] & \Fin_M \arrow[r, "M"] & \Set
    \end{tikzcd}
  \end{equation*}
\end{construction}

As MDet($M$) is based on the multiset monad it can be seen as a determinization
that remembers which states it passed and how many different ways it passed it.

We must again make sure that there is a canonical simulation between an automaton
and its determinization, as well as the fact that it recognizes the same
language.
Here the use of the local equivalence of categories and relative Kleisli-category
assures us of having a canonical simulation that is forward-backward, instead of
just forward.


\begin{lemma}
    Let $F: \mathcal{C}\to \SpanFin$ be an automaton,
    then there is a canonical forward-backward simulation from $F$ to
    $\CMamalg U F$ with pseudo naturality given by the $\eta$ of the local
    equivalence of categories.
    \begin{equation*}
        \begin{tikzcd}[ampersand replacement = \&]
            \mathcal{C} \arrow[rr, "F"] \& \& \SpanFin \arrow[rr, equal] \arrow[rd, "U"']
            \& {}
            \& \SpanFin \\
            \& {} \& {} \& \Fin_M \arrow[u, "{"\eta"}"', Rightarrow] \arrow[ru, "\CMamalg"']
            \&
        \end{tikzcd}
        \qquad
        \begin{tikzcd}[ampersand replacement = \&]
            A \arrow[rr, "\CMamalg_{(a,b)\in A\times B}{S(a,b)}"] \arrow[d, "\Id_{A}"']
            \& {} \& B \arrow[d, "\Id_B"] \\
            A \arrow[rr, "S"'] \arrow[rru, Rightarrow, "{\eta}", shorten =15pt]
            \& {} \& B
        \end{tikzcd}
    \end{equation*}
\end{lemma}

\begin{lemma}
    There is a natural transformation $\beta$, with the following naturality
    square (commuting in $\Set$)
    \begin{equation*}
        \begin{tikzcd}[ampersand replacement = \&]
            \& \Fin_M \arrow[rr, "i"] \arrow[rd, "U"']
            \& {}
            \& \Set_M \\
             \& {} \& \Set \arrow[u, "{\beta}"', Rightarrow] \arrow[ru, "\CMamalg"']
            \&
        \end{tikzcd}
        \qquad
        \begin{tikzcd}[ampersand replacement = \&]
            \mathbb{N}^A \arrow[rr, "\eta_{\mathbb{N}^B} f^*"] \arrow[d, "\Id_{\mathbb{N}^A}"']
            \& {} \& \mathbb{N}^{\mathbb{N}^B} \arrow[d, "\Id_{\mathbb{N}^B}^*"] \\
            \mathbb{N}^A \arrow[rr, "f^*"']
            \& {} \& \mathbb{N}^B
        \end{tikzcd}
    \end{equation*}
\end{lemma}

\begin{proposition}
    Let $F: \mathcal{C} \to \SpanFin$ be an automaton,
    then there is a canonical forward-backward simulation from $F$ to
    MDet($F$).
    \begin{equation*}
        \adjustbox{scale=0.88}{
        \begin{tikzcd}
            \mathcal{C} \arrow[rrr, "F"]
            & {} &
            & \SpanFin \arrow[rr, equal] \arrow[rd, "\text{U}"]
            & {} & \SpanFin \arrow[rr, "i"] & & \SpanSet \\
            & & {}
            & {}
            & \Fin_M \arrow[ru, "\CMamalg"'] \arrow[rr, "i"] \arrow[rd, "M"] \arrow[u, "{"\eta"}"', Rightarrow]
            & {} & \Set_M \arrow[ru, "i"'] & \\
            & {}
            & {}  &
            & & \Set \arrow[ru, "i"'] \arrow[u, "\beta"', Rightarrow]
            & &
            \end{tikzcd}}
    \end{equation*}
\end{proposition}

\subsection{Universal Property}

We again rely on constructing pasting inverses to obtain a universal property.
To construct these we rely on the fact that our construction is based on
a local equivalence of categories as well as a relative adjunction induced by
a relative monad.

\begin{lemma}\label{alpha}
    The natural transformation $\alpha$ is the identity when restricted to
    $\Fin_M$, as the $\epsilon$ from the local equivalence of categories
    \cref{localEQ} is the identity.
\end{lemma}

\begin{lemma}\label{beta}
    The natural transformation $\beta$ is the pasting inverse of the simulation
    $\eta$ from the relative adjunction
\end{lemma}

Using these lemmas we find the following corresponding universal property:

\begin{proposition}
    Let $G : \mathcal{C} \to \Set$ be some deterministic automaton, then any
    simulation from $F$ to $G$ factors through the canonical forward-backward
    simulation to MDet($F$) and a bisimulation between MDet($F$) and $G$.
\end{proposition}

\begin{proof}
    By \cref{alpha} and \cref{beta} we have that
    the following diagram:
    \begin{equation*}
        \adjustbox{scale=0.65}{
        \begin{tikzcd}
            \mathcal{C} \arrow[rdd, "G"', bend right] \arrow[rrr, "i F_p"]
            & {} &
            & \SpanSet\\
            & & \Fin_M \arrow[ru, "i \CMamalg"]\\
            & \Fin \arrow[ru, "i"]\arrow[uu, "\alpha"', Rightarrow]
            \end{tikzcd}}
    =
        \adjustbox{scale=0.65}{
        \begin{tikzcd}
            \mathcal{C} \arrow[rdd, "G"', bend right] \arrow[rrr, "F_p"]
            & {} &
            & \SpanFin \arrow[rr, equal] \arrow[rd, "\text{U}"]
            & {} & \SpanFin \arrow[rr, "i"] & & \SpanSet \\
            & & \Fin_M \arrow[ru, "\CMamalg"] \arrow[rd, "M"] \arrow[rr, equal]
            & {}
            & \Fin_M \arrow[ru, "\CMamalg"'] \arrow[rr, "i"] \arrow[rd, "M"] \arrow[u, "{"\eta"}"', Rightarrow]
            & {} & \Set_M \arrow[ru, "i"'] & \\
            & \Fin \arrow[ru, "i"] \arrow[rr, "i"] \arrow[uu, "\alpha"', Rightarrow]
            & {} \arrow[u, shorten <=0.8ex, "\eta"', Rightarrow] & \Set \arrow[rr, equal]
            & & \Set \arrow[ru, "i"'] \arrow[u, "\beta"', Rightarrow]
            & &
            \end{tikzcd}}
    \end{equation*}
    Consequently this gives a unique factorization of $\alpha$ through
    the universal simulation between $F$ and MDet($F$).
\end{proof}

\section{Conclusion}


In this paper, we studied determinizations for \Rel-automata and ULF-automata
from a categorical point of view. Determinization for \Rel-automata has
previously been investigated by Colcombet and Petri\c{s}an in \cite{colcombet},
as ULF-automata is a generalization of \Rel-automata, through the \SpanSet
intermediary, we consider how it might be generalized to \SpanSet-automata and
thus also to ULF-automata.

To do so, we gave a detailed presentation of the determinization procedure for \Rel, which
relies on the Kleisli adjunction between $\Set$ and \Rel, as
well as its correctness and universal property as previously described by
Colcombet and Petri\c{s}an.

Generalizing this construction to \SpanSet turned out to be non-trivial, as we
neither have an adjunction between \SpanSet and $\Set$, nor one between \SpanSet
and \Rel.
Moreover, as \SpanSet-automata corresponds to pseudofunctors, it forced us to
work in a 2-categorical setting and adapt definitions to it.
We have shown that it is possible to recover a fitting determinization procedure
by first turning a pseudo-functor $\mathcal{C} \to \SpanSet$ into $\mathcal{C} \to \Rel$
using the local adjunction and between \SpanSet and \Rel then using the
\Rel-determinization to get back a functor $\mathcal{C} \to \Set$.
To prove correctness and equip the procedure with a similar universal
property we used the local adjunction and the correctness and universal
property for \Rel-determinization.
The determinized automaton we obtained has subsets of states as states and a
transition function, which given a letter maps a subset to the set of
reachable states in the original automaton, analogously to
usual algorithmic determinization procedure.
We showed that our procedure does not always get whole powerset as the set of
states, but does as soon as we consider automata with a free monoid as codomain.

This determinization procedure has the particularity to identify the paths
taken, as this may not be wanted, we provided an alternative determinization
based on the multiset relative monad that retains them.
Compared to the one before, we no longer rely on the \Rel determinization
procedure as a middle step, as factoring through \Rel identifies paths.
Instead, we use the relative Kleisli category of the relative multiset monad $M$,
$\Fin_M$, as a middle step, and then postcompose by the functor from
$\Fin_M$ to $\Set$ provided by Kleisli construction of  $\Fin_M$.
Moreover, we have proven the correctness of this procedure and provided it with
a universal property, leveraging the relative monad and that there is a local
equivalence of categories between $\Fin_M$ and \SpanFin.
The multiset determinization enabled us to achieve path relevance, and to obtain
a stronger universal property that encompasses forward-backward simulation, as
we now get a canonical forward-backward simulation between an automaton $F$ and
its multiset determinization MDet($F$).

There are many potential ways to expand the work on categorical automata, one
potential possibility is to explore the multiset determinization from the point
of view of transducers where path relevance might be desired as we in
addition to the input also have an output to take into account.
Another possibility would be minimization, as this has been done for
\Rel-automata, we might consider how to generalize this to to our setting.

\section*{Acknowledgements}
The work for this paper was done during an summer internship under the
supervision of Noam Zeilberger, the author would thus like to thank him for the
introduction to the subject and many fruitful discussions.
The author is also grateful to Thomas Lamiaux for his advice on writing and
his comments on the many drafts of the paper.
Finally, the author would like
to thank the anonymous reviewers for their comments and suggestions.

\bibliographystyle{eptcs}
\bibliography{bib.bib}


\end{document}